\documentclass[a4paper,11pt,reqno]{amsart}

\usepackage[letterpaper, hmargin=1in, top=1in, bottom=1.2in, footskip=0.7in]{geometry}

\usepackage{rotating,pdflscape,color,amssymb,subfigure,psfrag,amsmath,eufrak,bbm,epsfig,amsthm, stmaryrd}
\definecolor{vdarkred}{rgb}{0.6,0,0.2}
\definecolor{vdarkblue}{rgb}{0,0.2,0.6}
\usepackage[pdftex, colorlinks, linkcolor=vdarkblue,citecolor=vdarkred]{hyperref}
\usepackage{a4wide,amscd}
\usepackage{tikz} 
 \usepackage[usenames,dvipsnames]{pstricks}
 \usepackage{pst-grad} 
 \usepackage{pst-plot} 

\usepackage[small]{caption}

\newcommand{\ii}{\mathrm{i}}
\newcommand{\me}{\mathrm{e}}

\newcommand{\cC}{\mathcal{C}}

\newcommand{\cS}{\mathcal{S}}

\newcommand{\bee}{\mathbf{e}}

\newcommand{\bu}{\mathbf{u}}

\newcommand{\wtD}{\widetilde{D}}

\newcommand{\wtX}{\widetilde{X}}

\newcommand{\ka}{\kappa} 

\newcommand{\si}{\sigma}
\newcommand{\lam}{\lambda}

\newcommand{\vfi}{\varphi}
\newcommand{\al}{\alpha}

\newcommand{\bet}{\beta}
\newcommand{\del}{\delta}

\newcommand{\ld}{\ldots}

\newcommand{\beg}{\begin}
\newcommand{\en}{\end}

\renewcommand{\Im}{\mathfrak{Im}}

\newcommand{\bgt}{\begin{itemize}}
\newcommand{\ent}{\end{itemize}}
\newcommand{\ite}{\item}
\newcommand{\eqre}{\eqref}
\newcommand{\re}{\ref}
\newcommand{\la}{\label}

\newcommand{\brem}{\begin{rmk}}
\newcommand{\erem}{\end{rmk}}
\newcommand{\blem}{\begin{lem}}
\newcommand{\elem}{\end{lem}}
\newcommand{\bcor}{\begin{cor}}
\newcommand{\ecor}{\end{cor}}
\newcommand{\bTh}{\begin{Th}}
\newcommand{\eTh}{\end{Th}}
\newtheorem{Claim}{Claim}[]
\newcommand{\bpropo}{\begin{propo}}
\newcommand{\epropo}{\end{propo}}

\newcommand{\sd}{\,\cdot\,}

\newcommand{\rfl}{\rfloor}
\newcommand{\lfl}{\lfloor}

\newcommand{\lan}{\langle}
\newcommand{\ran}{\rangle}

\newcommand{\op}{\operatorname}

\newcommand{\diag}{\operatorname{diag}}

\newcommand{\Tr}{\operatorname{Tr}}

\newcommand{\ud}{\mathrm{d}}

\newcommand{\ninf}{\underset{n\to\infty}{\longrightarrow}}

\newcommand{\E}{\op{\mathbb{E}}}

\newcommand{\R}{\mathbb{R}}
\newcommand{\C}{\mathbb{C}}

\newcommand{\p}{\mathbb{P}}

\newcommand{\pro}{probability }

\newcommand{\f}{\frac}
\newcommand{\ff}{\frac{1}}
\newcommand{\lf}{\left}
\newcommand{\ri}{\right}

\newcommand{\st}{such that }

\newcommand{\ti}{\times}

\newcommand{\mc}{\mathcal}

\newcommand{\eps}{\varepsilon}

\newcommand{\ovl}{\overline}

\newcommand{\bbm}{\begin{bmatrix}}
\newcommand{\ebm}{\end{bmatrix}}
\newcommand{\bes}{\begin{equation*}}
\newcommand{\ees}{\end{equation*}}
\newcommand{\be}{\begin{equation}}
\newcommand{\ee}{\end{equation}}
\newcommand{\beqy}{\begin{eqnarray}}
\newcommand{\eeqy}{\end{eqnarray}}
\newcommand{\beq}{\begin{eqnarray*}}
\newcommand{\eeq}{\end{eqnarray*}}
\newcommand{\one}{\mathbbm{1}}

\newcommand{\bpm}{\begin{pmatrix}}
\newcommand{\epm}{\end{pmatrix}}

\newcommand{\bpr}{\beg{proof}}
\newcommand{\epr}{\en{proof}}

\newcommand\numberthis{\addtocounter{equation}{1}\tag{\theequation}}

 \newcommand{\theo}[1]{Theorem \re{#1}}

\newtheorem{Th}{Theorem}

\newtheorem{propo}{Proposition}

\newtheorem{lem}{Lemma}

\newtheorem{cor}{Corollary}

\theoremstyle{definition}
\newtheorem{rmk}[]{Remark}

\long\def\symbolfootnote[#1]#2{\begingroup
\def\thefootnote{\fnsymbol{footnote}}\footnote[#1]{#2}\endgroup} 

\author[Florent Benaych-Georges]{Florent Benaych-Georges}\address{Florent Benaych-Georges, Universit\'e Paris Descartes,
45, rue des Saints-P\`eres
75270 Paris Cedex 06, France.}\email{florent.benaych-georges@parisdescartes.fr}

\author[Nathana\"el Enriquez]{Nathana\"el Enriquez}\address{Nathana\"el Enriquez, Universit\'e Paris-Sud,
Laboratoire Math\'ematiques d'Orsay, 91405 Orsay, France.} \email{nathanael.enriquez@u-psud.fr}

\author[Alk\'eos Micha\"{\i}l]{Alk\'eos Micha\"{\i}l}\address{Alk\'eos Micha\"{\i}l, Universit\'e Paris Descartes,
45, rue des Saints-P\`eres
75270 Paris Cedex 06, France.}\email{alkeos.michail@parisdescartes.fr}

\title{Eigenvectors of a matrix under random perturbation}

\date{\today}

\begin{document}

\maketitle
 
\begin{abstract}
In this text, based on elementary computations, we provide a perturbative expansion of the coordinates of the eigenvectors of a Hermitian matrix of large size perturbed by a random matrix with small operator norm whose entries in the eigenvector basis of the first one are independent, centered,  with a variance profile. This is done through a perturbative expansion of spectral measures associated to the state defined by a given vector.
\end{abstract}

\section{Introduction}

This paper is devoted to the study of the sensitivity of the eigenvectors of a given operator under small perturbations. 
In the previous paper \cite{BGEM} we studied the effect  of a perturbation on the spectrum of a diagonal matrix by a random matrix with small operator norm and whose entries in the eigenvector basis of the first one were independent, centered,  with a variance profile. We provided a perturbative expansion of the empirical spectral distribution, but   did not consider the deformation of the eigenvectors basis with respect to the canonical basis. In the present paper,   to complete this first study, we  deal with the spectral measure of our matrix associated to the   state defined by a given vector.

To define this measure, let us introduce some notations. We consider a real diagonal matrix $D_n =\diag(\lam_{1}, \ld, \lam_{n})$ (the eigenvalue $\lam_{i}$ implicitly depends on $n$),  as well as a Hermitian random matrix $$X_n=\ff{\sqrt{n}}\lf[x^n_{i,j}\ri]_{1\le i,j\le n}$$ \st the $x_{ij}$ are independent (up to the symmetry), centered,  with a variance profile. The normalizing factor  $ n^{-1/2}$ and our hypotheses below ensure that the operator norm of $X_n$ is of order one.
We then define,   for   $\eps>0$,
$$D_n^\eps:=D_n+ \eps X_n.$$

If the perturbing matrix belongs to the GOE or GUE, then its law is invariant under this change of basis, hence all the results of this paper apply to any self-adjoint matrix $D_n$.

In contrast with \cite{BGEM}, where we studied the empirical spectral measure $\mu_n^\eps$ of the matrix $D_n^\eps$,    we consider here   the spectral measure $\mu_{n, \bee_i}^\eps$  of $D_n^\eps$ over a vector $\bee_i$ of the canonical basis, defined\newline\newline through an orthonormal eigenbasis $(\bu_j^\eps)_{j\in\{1, \dots, n \} }$ of $D_n^\eps$  and the related eigenvalues $(\lambda_j^\eps)_{j\in\{1, \dots, n \} }$  by
\[ \mu_{n,\bee_i}^\eps:=\sum_{j=1}^n |\lan \bu_j^\eps, \bee_i\ran|^2  \del_{\lam_j^\eps}.\] 

The interest of these measures is  that they give information on the eigenvector basis of $D_n^\eps$, while being tractable since they satisfy, for any test function $\vfi$, the   key identity
\be\label{keyIdentity}
 \int \vfi(x) \ \ud \mu_{n,\bee_i}^\eps(x) =
\sum_{j=1}^n |\lan \bu_j^\eps, \bee_i\ran|^2  \vfi(\lam_j^\eps) = \left( \vfi\left( D_n^\eps \right) \right)_{i,i}.
\ee

Our main result, Theorem \re{theoMain0917},  gives a perturbative expansion of  $\mu_{n,\bee_i}^\eps$. More precisely,  using   a resolvent expansion and the Helffer-Sj\"ostrand formula, we give   an asymptotic expansion of $$ \int_\R\vfi(t)\ud \mu_{n,\bee_i}^\eps(t)$$ for any $\cC^5$ test function $\vfi$. From that, we deduce \theo{theoMain0918} which establishes the convergence of the average of the square  of coordinates of a mesoscopic sequence of consecutive eigenvectors.

 It would be indeed tempting to generalize this analysis to the non-diagonal entries of the matrix $\varphi(D_n^\eps)$. For $1 \leq l,k \leq n$ the entry $\varphi(D_n^\eps)_{k,l}$ would give access to the measure $ \sum_{j=1}^{n} \lan \bu_j^\eps, \bee_l\ran \lan \bu_j^\eps, \bee_k\ran \delta_{\lam_j^\eps} $. A result on its asymptotic behavior  as the one we have for $k=l$ can not lead to more information than the mere intensity of $\lan\bu_j^\eps, \bee_l\ran$ and $\lan\bu_j^\eps, \bee_k\ran$ \textit{separately} in terms of the distances $|j-l|$ and $|j-k|$. Information on correlations is beyond what we can get with our method.

Some  other works, on models closed to our one or contained in it,  are   devoted  to the sensitivity to perturbations of the eigenvectors. Some of them, as \cite{ORVK0,ORVK,PVSSW,Zhong2017}, provide bounds on the deviations of these eigenvectors under perturbation,  while some other, as \cite{allez2012eigenvector,allez2014eigenvector,AllezBunBouchaud,Benigni},   provide explicit perturbative expansions.  
This is what we do here,
our Theorem \ref{theoMain0918}  shows that the   overlaps   $|\lan \bu_j^\eps, \bee_i\ran|^2$ have order $\eps^2(\lam_{j}-\lam_{i})^{-2} n^{-1}$. We cannot prove it for all indices $i,j$ individually but only  in average  over some mesoscopic windows. The size of window we have to take is larger than $n^{11/12}$ which is certainly not optimal as suggested by the recent work of Benigni \cite{Benigni} which has very refined  and non perturbative in $\eps$ results in the special case when $X_n$ is a Wigner matrix. He makes use, among others, of the sophisticated method of Bourgade and Yau called the eigenvector moment flow \cite{BY}.
In addition to the fact that it only relies on  short and elementary computations, one of the interests of the present paper is  to consider rather general perturbations, since we do not suppose that all entries of $X_n$ have the same variance nor that they are Gaussian. Another interest is to provide, with the functional $\Xi_s(\vfi) $ from \eqre{defXi_s} and \eqre{3101710}, an   expression for the first order expansion of the measure $\mu_{n,\bee_i}^\eps$ from \eqre{keyIdentity}, which, up to our knowledge, did not appear so far.

The   paper is organized as follows. Statement of Theorem \re{theoMain0917} and comments are given in Section \ref{sectionMainResutlEigenvectors}, whereas its proof is given in  Section \ref{sectionProofEigenvectors}. Section \ref{sectionHeuristic}
is devoted to the consequence  of Theorem \re{theoMain0917} on   the eigenvectors,  namely to  Theorem \ref{theoMain0918}, some comments on this result and some figures.   Theorem \re{theoMain0918} is proved in Section \ref{sectionHeuristicproof}.

{\bf Notations.} For $u=u_n$, $v=v_n$ some   sequences, $u\ll v$  
means that $u_n/v_n$ tends to $0$.\\
For a given sequence $u_n$ we denote by $O_{L^2}(u_n)$ any sequence $U_n$ of random variables whose $L^2$ norm $\E(U_n^2)^{1/2}$ is uniformly bounded by $C u_n$ for some $C>0$.\\
Finally, we denote by $\xrightarrow[n\to\infty]{P}$ the convergence in probability for sequences of random variables.

\section{Main result}
\label{sectionMainResutlEigenvectors}
  
   We consider a real diagonal matrix $D_n =\diag(\lam_{1}, \ld, \lam_{n})$ (the eigenvalue $\lam_{i}$ implicitly depend on $n$),  as well as a Hermitian random matrix $$X_n=\ff{\sqrt{n}}\lf[x^n_{i,j}\ri]_{1\le i,j\le n}$$ and  define,   for   $\eps=\eps_n>0$,
$$D_n^\eps:=D_n+ \eps X_n.$$

We make the following hypotheses:
 \begin{enumerate}\item[(a)]  the entries $x^n_{i,j}$ of $\sqrt{n}X_n$ are independent (up to symmetry) random variables, centered, with variance denoted  by $\si_n^2(i,j)$, such that $\E   |x^n_{i,j}|^{6} $ is bounded uniformly on $n,i,j$,\\
  \item[(b)] there are two bounded real functions, $f$ and $\si$, defined respectively on $[0,1]$   and $[0,1]^2$ such that,  denoting $\lam_i$ by $\lam_{n,i}$ to emphasize the implicit dependence in $n$,  the error bound\beqy\la{defetan71017}  \eta_n&:=&\sup_{x\in [0,1]}|\lam_{n,\lfl nx\rfl}- f(x)|  +\sup_{(x, y)\in [0,1]^2}|
  \si_n^2(\lfl nx\rfl,\lfl ny\rfl)- \si^2(x,y)|\eeqy
  satisfies $$\eta_n\ninf 0.$$
  \end{enumerate}
  Let us now make some   assumptions on the limiting functions  $\si$ and $f$: \begin{enumerate}\item[(c)]   the push-forward   of the uniform measure on $[0,1]$ by the function $f$ has a  density $\rho$ with respect to the Lebesgue measure on $\R$ and a   compact support denoted by $\cS$,\\
 \item[(d)] the variances of the entries of $X_n$ essentially depend on the eigenspaces of $D_n$, namely, there exists a symmetric    function $\tau(\sd,\sd )$ on $\R^2$ \st for all $x\ne y$, $\si^2(x,y)=\tau(f(x),f(y))$.
\end{enumerate}

\begin{rmk}
We refer the reader to the end of Section \ref{sectionHeuristic} for matrix models satisfying these hypotheses.
\end{rmk}

\begin{rmk} The part of assumption $(a)$ concerning the existence of the sixth moment of $x_{i,j}$ is due to our aim at giving a Taylor type expansion of the Stieltjes transform of the spectral measure. In this respect it is very likely to be optimal.
Assumption $(c)$ prevents us considering the case when $D_n$ is a scalar matrix since its limiting empirical measure has no density.
\end{rmk}

\begin{rmk}
We cannot generalize our result to the case when $D_n$ is Hermitian non diagonal, unless the perturbations belong to GOE or GUE. The reason is mainly due to assumption $(a)$ of independence of the entries of $X_n$. It seems challenging to study the more general problem assuming the existence of a limiting correlation profile. Looking carefully at the proof of Claim 2 shows that correlations of order $\tfrac{1}{n}$ do not change the magnitude of our error terms and that we can maintain a statement as long as the correlations are of order $o(1)$.
\end{rmk}

Let $ \mu_{n,\bee_i}^\eps$
denote the \pro measure defined, for any test function $\vfi$, by 
\be\la{defmui1}\int\vfi(t)\ud \mu_{n,\bee_i}^\eps(t):=(\vfi(D_n^\eps))_{ii}.\ee 
One can  equivalently define $\mu_{n,\bee_i}^\eps$ by 
\be\la{defmui2}\mu_{n,\bee_i}^\eps:=\sum_{j=1}^n |\lan \bu_j^\eps, \bee_i\ran|^2\del_{\lam_j^\eps},\ee 
where $\bee_i$ denotes the $i$-th vector of the canonical basis, the $\lam_j^\eps$'s denote the eigenvalues of $D_n^\eps$ and the $\bu_j^\eps$'s denote the associated eigenvectors.

 We now   introduce a functional which is central in the statement of our result. This functional admits another expression, given in Proposition  \re{propo:rewritingFunctional} below.
  
  Let, for $s\in \R$ and $\vfi:\R\to\C$ a $\cC^2$ function,
  \be\la{defXi_s} \Xi_s(\vfi):=  \int_{\mathbb{R}}\tau(s,t)\rho(t)\f{(\vfi(t)-\vfi(s)-(t-s)\vfi'(s))}{(t-s)^2} \ud t\ee
  
  \bTh\la{theoMain0917} Let us suppose that $\eps=\eps_n\ll n^{-\frac12}$.  Let  $\vfi:\R\to\C$ be a compactly supported $\cC^7$ function. For $x\in[0,1]$, set $i=i(n,x)=\lfloor nx \rfloor$. Then  we have $$ \int_\R\vfi(t)\ud \mu_{n,\bee_i}^\eps(t)=\vfi\left(\lam_{i}+\f{\eps}{\sqrt n} x_{ii}\right)+\eps^2\Xi_{f(x)}(\vfi)+\eps^2  O_{L^2}\left( \|\vfi^{(7)}\|_\infty \ \left(\eta_n + \frac1{\sqrt{n}} + \eps \sqrt{n} \right)\right)$$ 
     for $\eta_n$ as in \eqre{defetan71017}.\eTh
 
 \brem[Leading order transition]\la{rmkLO} Note that for any $\cC^2$ test function  $\vfi$, 
$$\vfi\lf(\lam_{i}+\f{\eps}{\sqrt n} x_{ii}\ri)=\vfi(\lam_{i})+\f{\eps}{\sqrt n} x_{ii}\vfi'(\lam_{i})+O_{L^2}\lf( \f{\eps^2}{n} \|\vfi''\|_\infty\ri). $$
Thus the previous theorem allows to expand the measure $\mu_{n,\bee_i}^\eps$ around $\del_{\lam_i}$ as follows. With the notations and the hypothesis of the theorem, 
\beqy\nonumber  \int\vfi(t)\ud \mu_{n,\bee_i}^\eps(t)&=&\vfi(\lam_{i})+\f{\eps}{\sqrt n} x_{ii}\vfi'(\lam_{i})+\eps^2\Xi_s(\vfi)\\ &&\qquad + O_{L^2}\lf(\eps^2\|\vfi^{(5)}\|_\infty(\eps n^\frac12+n^{-\frac12}+\eta_n)+ \f{\eps^2}{n} \|\vfi''\|_\infty\ri).\la{eqRem1260118}\eeqy

If $\vfi'(\lam_{i}) \neq 0$ then the assumption $\eps \ll n^{-\frac{1}{2}}$ implies that the term $\f{\eps}{\sqrt n} x_{ii}\vfi'(\lam_{i})$  prevails over the term $ \eps^2\Xi_s(\vfi) $ but in the following, we will apply Theorem \ref{theoMain0917} to test functions whose support avoids $\lam_{i}$, so that $\eps^2\Xi_s(\vfi)$ will be the dominant term of the expansion.

\end{rmk}

\begin{rmk}\label{remarkSimpleXi}
Strikingly, the image of a function $\vfi$ by the operator $\Xi_{f(x)}$ is not changed if one adds an affine function to $\vfi$. 
This can be understood because the measure
$\mu_{n,\bee_i}^\eps  - \delta_{\lambda_i}$
is of null mass and with first moment of order $o(\eps^2)$ since by \eqref{keyIdentity},
\[ \int_\R x\ \ud\left(\mu_{n,\bee_i}^\eps  - \delta_{\lam_{i}}\right) = \left(D_n^\eps\right)_{ii} - \lam_{i}= \frac{\eps}{\sqrt{n}} x_{ii} = o(\eps^2). \]
Note that when both $\vfi(f(x))$ and $\vfi'(f(x))$ are null, the function $\Xi_{f(x)}(\vfi)$ boils down to the integral
\beqy\label{SimpleXi}
\Xi_{f(x)}(\vfi) =  \int_{\mathbb{R}}\tau(f(x),t)\rho(t)\f{\vfi(t)}{(t-f(x))^2} \ud t.
\eeqy
We will use this fact in Section \ref{sectionHeuristic} for test functions $\vfi$ whose support does not contain $f(x)$.
\end{rmk}

 \bpropo\la{propo:rewritingFunctional}
Let us define, for any $s\in\R$, the function $\zeta_s$ defined on $\mathbb{R}$ by
\begin{equation}\label{zeta_s}
\zeta_s(y) := \int_{1}^{+\infty}\f{r-1}{r^2}\tau(s,s+r(y-s)) \rho(s+r(y-s))\ud r.
\end{equation}
Then for any $\cC^2$ function $\vfi$ and any $s\in \R$, the functional $\Xi_s$ defined at \eqre{defXi_s}  rewrites
  \beqy \la{3101710}\Xi_s(\vfi)
&=&\int_{\R}\vfi''(y) \zeta_s(y) \ud y.
\eeqy  
  \epropo
  
  \bpr Taylor's formula yields $$\vfi(t)-\vfi(s)-(t-s)\vfi'(s)=\int_{s}^t\vfi''(x)(t-x)\ud x=(t-s)^2\int_{u=0}^1\vfi''(s+u(t-s))(1-u)\ud u.$$
  Hence,
  \[\Xi_s(\vfi)  = \int_{t\in \R}\int_{u=0}^1\vfi''(s+u(t-s))(1-u)\ud u\ \tau(s,t) \rho(t)\ud t\]
  We now perform the change of variable $(r,y) = \Psi_s(u,t)$ with 
  $$\Psi_s:(u,t)\in(0,1)\ti\R\mapsto (r,y)=\left(\frac{1}{u},u(t-s)+s\right)\in (1,\infty)\ti\R$$
   which gives the result
  \epr

  \section{Consequence for the eigenvectors}\label{sectionHeuristic}

The purpose of this section is to use the previous results to obtain information on the projection of the eigenvectors on the canonical basis (via moving averages of course, as  seeking to obtain a result about eigenvectors one by one would be unrealistic at this level of generality).

\bTh\la{theoMain0918} 
For all sequences $\al_n$ converging to zero and satisfying $\al_n^8 \gg \max\left\{n^{\frac12}\eps , \eta_n, n^{-\frac12}\right\}$, for all $x, x_0\in[0,1]$ with $x\neq x_0$, the following convergence in probability holds,
$$ \frac{n\eps^{-2}}{\op{Card} \{ j\,:\, |\lam_j^\eps - f(x)| < \al_n  \}} \sum_{\left\{ j\,:\, |\lam_j^\eps - f(x)| < \al_n \right\}}  |\lan\bu_{j}^\eps,\bee_{\lfloor nx_0 \rfloor}\ran|^2 \xrightarrow {{ {P}}} \frac{\tau(f(x_0), f(x))}{\left( f(x) - f(x_0) \right)^2}. $$
\eTh

\begin{rmk}
This is a \emph{local} result since the  window where we take our average contains  $o(n)$ eigenvectors. However, this $o(n)$  is at least $n^{15/16}$, which is for sure not optimal, as suggested by the recent work of Benigni \cite{Benigni} who gets a very refined result in the special case where the perturbating matrix is Wigner (which implies,  among other, that $\tau\equiv 1$). He proves actually that the components of the eigenvectors are asymptotically independent and normal and gets therefore the convergence in probability for any size of window converging to infinity.
\end{rmk}

We present now two simulations (displayed in Figures \re{figure1Eigenvectors} and \re{figure2Eigenvectors}) which show a   good matching with this theoretical prediction.  First we consider the case where the deterministic matrix $D_n$ is perturbed by a Gaussian Wigner matrix, $X_n$. More precisely, we take for $D_n$ the diagonal matrix with $\frac{i}{n}$ as $i^{th}$ entry, so that $f(x)=x$ and the density $\rho$ is equal $x\mapsto \one_{[0,1]}(x)$. The entries of the perturbating matrix $X_n$ are all Gaussian and independent with variance one. Then, we consider the case where the same matrix $D_n$ is perturbed  by a band matrix. In other words, we consider now that $\sigma(x,y)=\one_{|x-y|\leq \ell}$, where $\ell\in[0,1]$ is the relative width of the band. Note that in this second example, even though there is absolutely no deterministic reason why $\lan\bu_{\lfloor n y \rfloor}^\eps,\bee_{\lfloor nx \rfloor}\ran$ would vanish when  $|y-x|>\ell$, we see that at first order, it is actually almost zero (Figure  \re{figure2Eigenvectors}). This is related to the question of the localization of the eigenvectors of band random matrices  (see e.g. \cite{DjalilAnderson,EKYY,EK2011CMP,ErdKnow,EYY,EYY2,FM,Schenker}). \begin{figure}[!h]
\centering
\includegraphics[scale=0.65]{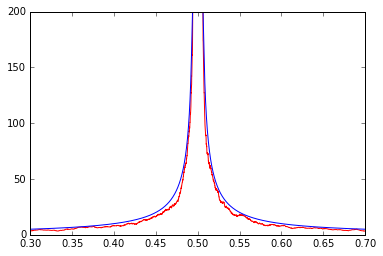}
\caption{ \textbf{Uniform measure perturbation by a Wigner matrix.} The red curve represents  a moving average of the function $t \in [0,1]\longmapsto\eps^{-2} n |\lan\bu_{\lfloor n f^{-1}(t) \rfloor}^\eps,\bee_{\lfloor nx_0 \rfloor}\ran|^2,$ over a  window of length  $\frac1{\sqrt n}$. The blue curve represents our theoretical prediction $t\longmapsto  |t-x_0|^{-2}$.
Here $n=10^4$, $\varepsilon = n^{-0.7}$ and $x_0=\tfrac{1}{2}$.}
\label{figure1Eigenvectors}
\end{figure}

\begin{figure}[!h]
\centering
\includegraphics[scale=0.65]{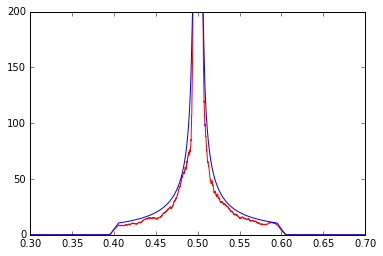}
\caption{ \textbf{Uniform measure perturbation by a band matrix.} The red curve represents  a moving average of the function $t \in [0,1]\longmapsto\eps^{-2} n |\lan\bu_{\lfloor n f^{-1}(t) \rfloor}^\eps,\bee_{\lfloor nx_0 \rfloor}\ran|^2,$ over a  window of length  $\frac1{\sqrt n}$. The blue curve represents our theoretical prediction $t\longmapsto \one_{|t-f(x_0)|\leq \ell} |t-x_0|^{-2}$.
Here $n=10^4$, $\ell=0.1$, $\varepsilon = n^{-0.7}$ and $x_0=\tfrac{1}{2}$.}
\label{figure2Eigenvectors}
\end{figure}

\section{Proof of \theo{theoMain0917}}
\label{sectionProofEigenvectors}

The proof is divided into two parts.
We shall first prove a convergence result for test functions $\varphi$ of the type $\vfi_z:=\ff{z-x}$. This is the purpose of Subsection \ref{subsec1311}. It will be obtained by writing an expansion of the resolvent of $D_n^\eps$. 

Once we have proved that such a convergence holds for the resolvent of $D_n^\eps$, we will be able to extend it to the class of compactly supported $\mathcal{C}^7$ functions on $\R$, by using the Helffer-Sj\"ostrand formula (see \cite{Helffer1989} or \cite{NotesAnttiFlo}) which expresses a regular function $\varphi$ on $\R$ as an integral against functions  $\vfi_z$ of the previous type. This is done in Subsection \ref{subsec1312}.

\subsection{Stieltjes transform}\label{subsec1311}
Let us introduce the Banach space $\mc{C}^2_{\op{b}}$ of bounded  $\mc{C}^2$ functions on $\R$ with  bounded first and second derivatives, endowed with the  norm $\|\vfi\|_{\mc{C}^2_{\op{b}}}:=\|\vfi\|_\infty+\|\vfi'\|_\infty+\|\vfi''\|_\infty$.

On this space, let us define, for $x\in[0,1]$ and $i = \lfloor  nx \rfloor$, the random continuous  linear form
\bes\label{30nov1554} \Pi_n(\vfi):= \eps^{-2}\lf(\int\vfi(t)\ud \mu_{n,\bee_i}^\eps(t)-\vfi(\lam_i+\f{\eps}{\sqrt n} x_{ii})\ri)- \Xi_{f(x)}(\vfi). \ees 

\begin{lem}\label{lemmaPi_n}
There exists a constant $C>0$ \st for all $z\in\C\setminus\R$,
\begin{align}
\E[|\Pi_n(\varphi_z)|^2] \; &\leq \; C \lf( \frac{\eta_n^2 + \frac{1}{n}}{|\Im(z)|^6}  + \frac{n\eps^2}{ |\Im(z)|^8} + \frac{\eps^4}{n^2 |\Im(z)|^{10}}  + \frac{\eps^6}{n^{3} |\Im(z)|^{12}} \ri). \label{1504091216}
\end{align}
\end{lem}

\begin{rmk}
This result implies that $\forall z\in\C\setminus\R, \ \Pi_n(\varphi_z) \xrightarrow[n\to\infty]{P} 0.$
\end{rmk}

Let us prove the above lemma. We denote, for short, $x^n_{i,j}$ by $x_{ij}$
and introduce the diagonal matrix
\beqy \label{wDeps}
\wtD_n^\eps:=\diag\left( \left( \widetilde{\lam}_n^\eps(i) :=  \lam_i+\f{\eps}{\sqrt n}x_{ii}\right)_{ i=1, \ld, n}\right)
\eeqy
which is the diagonal part of the matrix $D_n^\eps$. Note that with this notation and by using identity \eqref{keyIdentity}, the quantity we are interested in can be written:
\beqy
\label{Pi_nInterest}
\Pi_n(\vfi_z) = \eps^{-2} \lf((z-D_n^\eps)^{-1}-(z-\wtD_n^\eps)^{-1}\ri)_{ii} - \Xi_{f(x)}(\vfi_z).
\eeqy

To deal with this quantity we introduce the null diagonal matrix 
$$\wtX_n:= \eps^{-1}(D_n^\eps - \wtD_n^\eps) = X_n-n^{-1/2}\diag((x_{ii})_{i=1, \ld, n})$$
obtained by vanishing the diagonal of the matrix $X$.

A perturbative expansion of the resolvent of $D_n^\eps = \wtD_n^\eps + \eps\wtX_n$ yields
\beqy
\label{diagonalExpansion}
 (z-D_n^\eps)^{-1}-(z-\wtD_n^\eps)^{-1}&=&\eps(z-\wtD_n^\eps)^{-1}\wtX_n (z-\wtD_n^\eps)^{-1}
\\ &&+\eps^2 (z-\wtD_n^\eps)^{-1} \wtX_n(z-\wtD_n^\eps)^{-1}\wtX_n(z-\wtD_n^\eps)^{-1}\nonumber \\ &&
+\eps^3 (z-\wtD_n^\eps)^{-1} \wtX_n(z-\wtD_n^\eps)^{-1}\wtX_n(z-\wtD_n^\eps)^{-1}\wtX_n(z- D_n^\eps)^{-1}.\nonumber \eeqy

We now want to analyze the corresponding expansion of 
$ 
\lf((z-D_n^\eps)^{-1}-(z-\wtD_n^\eps)^{-1}\ri)_{ii}.
$

\begin{Claim}
For all $i \in \llbracket 1,n \rrbracket$,
$\lf((z-\wtD_n^\eps)^{-1}\wtX_n (z-\wtD_n^\eps)^{-1}\ri)_{ii} = 0$.
\end{Claim}

\bpr
This comes from the fact that the matrix $\wtX_n$ has a null diagonal.
\epr

\begin{Claim}
If, for all $i \in \llbracket 1,n \rrbracket$, we denote
 $$B_n(z,i):= \left(  (z-\wtD_n^\eps)^{-1} \wtX_n(z-\wtD_n^\eps)^{-1}\wtX_n(z-\wtD_n^\eps)^{-1}  \right)_{ii},$$
then, for all $x\in[0,1]$, 
$$ B_n(z,\lfloor nx \rfloor) - \Xi_{f(x)}(\vfi_z) = O_{L_2}\lf( \f{\eta_n + \ff{\sqrt{n}}}{|\Im(z)|^3}  + \frac{\eps}{\sqrt{n} |\Im(z)|^4} + \frac{\eps^2}{n |\Im(z)|^5}  + \frac{\eps^3}{n^{\frac{3}{2}} |\Im(z)|^6} \ri).$$
\end{Claim}

\bpr
With  the notations of \eqref{wDeps}, the term $B_n(z,i)$ writes
$$B_n(z,i) = \ff{n}\sum_{j=1}^n\f{|x_{ij}|^2}{\lf(z-\widetilde{\lam}_n^\eps(i)\ri)^2 \lf(z-\widetilde{\lam}_n^\eps(j)\ri)},$$ 
and, for $x\in [0,1]$,
\beqy
\Xi_{f(x)}(\vfi_z) &=& 
\int_{t\in\R} \frac{\tau(f(x),t) \rho(t)}{(t-f(x))^2} \left( \frac1{z-t} - \frac1{z-f(x)} - \frac{t-f(x)}{(z-f(x))^2} \right) \ud t \nonumber
\\
&=& \int_{t\in \R}\f{\tau(f(x),t)}{(z-f(x))^2(z-t)}\rho(t)\ud t  \la{310171} 
\\
&=& \int_{y\in [0,1]}\f{\si^2(x,y)}{(z-f(x))^2(z-f(y))} \ud y. \la{310172}
\eeqy
The difference of these quantities writes,
\begin{align*}
& B_n(z, \lfloor nx \rfloor) - \Xi_{f(x)}(\vfi_z) 
\\
& = \ \ff{n}\sum_{j=1}^n \f{|x_{\lfloor nx \rfloor,j}|^2-\si_n^2\lf(\lfloor nx \rfloor,j\ri)}{\lf(z - \widetilde{\lam}_n^\eps(\lfloor nx \rfloor)\ri)^2 \lf(z-  \widetilde{\lam}_n^\eps(j) \ri)}
\\
& \quad + \frac1n \sum_{j=1}^n  \frac{\si_n^2\lf(\lfloor nx \rfloor,j\ri)}{\lf(z - \widetilde{\lam}_n^\eps(\lfloor nx \rfloor) \ri)^2 \lf(z- \widetilde{\lam}_n^\eps(j) \ri)} - \frac{\si_n^2\lf(\lfloor nx \rfloor,j\ri)}{\lf(z-\lam_{\lfloor nx\rfloor}\ri)^2 \lf(z- \lam_j \ri)}
\\
& \quad + \int_{y\in [0,1]}\f{\si_n^2\lf( \lfloor nx \rfloor,\lfl ny\rfl\ri)}{(z-\lam_{\lfloor nx \rfloor})^2  (z-\lam_{\lfloor ny\rfloor})} \ud y -\int_{y\in [0,1]}\f{\si^2(x,y)}{(z-f(x))^2(z-f(y))} \ud y.
\end{align*}

Observe that the second-to-last integral coincides with the discrete sum $ \frac1n \sum_{j=1}^n  \frac{\si_n^2\lf(\lfloor nx \rfloor,j\ri)}{\lf(z-\lam_{\lfloor nx\rfloor}\ri)^2 \lf(z- \lam_j \ri)}$ since it concerns step functions.

Using the key assumption about the \textit{independence} of the variables $(x_{i,j})$, the $L_2$ norm of the first line of the right hand side of the previous equality writes
 $$ \lf\|\ff{n}\sum_{j=1}^n \f{|x_{ij}|^2-\si_n^2(i,j)}{(z-\lam_i)^2(z-\lam_j)}\ri\|_{L^2} =\ff{n} \lf( \sum_{j=1}^n \f{\E[(|x_{ij}|^2-\si_n^2(i,j))^2]}{|z-\lam_i|^4|z-\lam_j|^2} \ri)^{\frac{1}{2}} = O\lf(\ff{\sqrt{n} \ |\Im z|^{3}} \ri),$$
to analyze the $L^2$ norm of the second line, we write
\begin{align*}
& \frac1{\lf(z - \widetilde{\lam}_n^\eps(\lfloor nx \rfloor) \ri)^2 \lf(z- \widetilde{\lam}_n^\eps(j) \ri)} - \frac1{\lf(z-\lam_{\lfloor nx\rfloor}\ri)^2 \lf(z- \lam_j \ri)}
\\
= & \frac{\eps}{\sqrt{n}} \  \frac{ 2\lf(z-\lam_{\lfloor nx\rfloor}\ri)  \lf(z-\lam_{j}\ri)   x_{\lfloor nx \rfloor, \lfloor nx \rfloor} + \lf(z-\lam_{\lfloor nx\rfloor}\ri)^2 x_{j,j}  }{\lf(z - \widetilde{\lam}_n^\eps(\lfloor nx \rfloor) \ri)^2 \lf(z- \widetilde{\lam}_n^\eps(j) \ri)  \lf(z-\lam_{\lfloor nx\rfloor}\ri)^2 \lf(z- \lam_j \ri)}
\\
- & \frac{\eps^2}{n} \ \frac{   x_{\lfloor nx \rfloor, \lfloor nx \rfloor}^2 \lf(z-\lam_{j}\ri)  + 2 \lf(z-\lam_{\lfloor nx\rfloor}\ri)  x_{\lfloor nx \rfloor, \lfloor nx \rfloor} x_{j,j} }{\lf(z - \widetilde{\lam}_n^\eps(\lfloor nx \rfloor) \ri)^2 \lf(z- \widetilde{\lam}_n^\eps(j) \ri)  \lf(z-\lam_{\lfloor nx\rfloor}\ri)^2 \lf(z- \lam_j \ri)}
\\
+ & \frac{\eps^3}{n^{\frac{3}{2}}} \ \frac{  x_{\lfloor nx \rfloor, \lfloor nx \rfloor}^2  x_{j,j}    
}{\lf(z - \widetilde{\lam}_n^\eps(\lfloor nx \rfloor) \ri)^2 \lf(z- \widetilde{\lam}_n^\eps(j) \ri)  \lf(z-\lam_{\lfloor nx\rfloor}\ri)^2 \lf(z- \lam_j \ri)},
\end{align*}
hence,
\begin{align*}
&\left| \frac1{\lf(z - \widetilde{\lam}_n^\eps(\lfloor nx \rfloor) \ri)^2 \lf(z- \widetilde{\lam}_n^\eps(j) \ri)} - \frac1{\lf(z-\lam_{\lfloor nx\rfloor}\ri)^2 \lf(z- \lam_j \ri)} \right|
\\
\leq &
\frac{\eps}{\sqrt{n}} \frac{ 2| x_{\lfloor nx \rfloor, \lfloor nx \rfloor}| + |x_{j,j}|}{ |\Im(z)|^4}
+
\frac{\eps^2}{n} \frac{ 2| x_{\lfloor nx \rfloor, \lfloor nx \rfloor}|^2 + |x_{j,j}|^2}{ |\Im(z)|^5}
+ \frac{\eps^3}{n^{\frac{3}{2}}} \frac{ | x_{\lfloor nx \rfloor, \lfloor nx \rfloor}|^2  |x_{j,j}|}{ |\Im(z)|^6}.
\end{align*}
Therefore,\\
\\
$
\left\| \frac1n \sum_{j=1}^n  \frac{\si_n^2\lf(\lfloor nx \rfloor,j\ri)}{\lf(z - \widetilde{\lam}_n^\eps(\lfloor nx \rfloor) \ri)^2 \lf(z- \widetilde{\lam}_n^\eps(j) \ri)} - \frac{\si_n^2\lf(\lfloor nx \rfloor,j\ri)}{\lf(z-\lam_{\lfloor nx\rfloor}\ri)^2 \lf(z- \lam_j \ri)}
\right\|_{L^2} = O\left( \frac{\eps}{\sqrt{n} |\Im(z)|^4} + \frac{\eps^2}{n |\Im(z)|^5}  + \frac{\eps^3}{n^{\frac{3}{2}} |\Im(z)|^6}  \right).
$
and, finally, from assumption (b), the third line is $O(\eta_n|\Im z|^{-3}).$
\epr

\begin{Claim}\label{Claim3}
For all $i\in\llbracket 1, n \rrbracket$,
$$ \lf( (z-\wtD_n^\eps)^{-1} \wtX(z-\wtD_n^\eps)^{-1}\wtX(z-\wtD_n^\eps)^{-1}\wtX(z-D_n^\eps)^{-1} \ri)_{ii} =O_{L^2}(n^{\frac12} |\Im z|^{-4})$$
\end{Claim}

\bpr

If we denote $E_{i,i}$ the matrix with null entries everywhere except in position $(i,i)$ where the entry is equal to $1$, the Cauchy-Schwarz inequality yields

\begin{align*}
&\E\left[  \left| \left((z-\wtD_n^\eps)^{-1} \wtX_n (z-\wtD_n^\eps)^{-1}\wtX_n (z-\wtD_n^\eps)^{-1}\wtX_n (z- D_n^\eps)^{-1} \ri)_{ii}\right|^2 \right]^{\ff{2}}
\\
&= \ \E\lf[\lf|\Tr \lf( E_{i,i} \lf( \ff{z-\wtD_n^\eps} \wtX_n   \ff{z-\wtD_n^\eps} \wtX_n    \ff{z-\wtD_n^\eps}\wtX_n   \ff{z-D_n^\eps} \ri) \ri) \ri|^2\ri]^{\ff{2}}
\\
&\leq \ \E\lf[ \Tr \lf|E_{i,i}\lf(\ff{z-\wtD_n^\eps}\wtX_n \ri)^{3}\ri|^{2} \times \Tr \lf| \ff{z-D_n^\eps} \ri|^2 \ri]^{\ff{2}} \numberthis \label{TracesProduct}
\end{align*}

Let us now observe that since the spectra of $\wtD_n^\eps$ and $D_n^\eps$ are real, the moduli of the entries of $(z-\wtD_n^\eps)^{-1}$ and of  $(z-D_n^\eps)^{-1}$ are smaller than $|\Im(z)|^{-1}$, which implies that,
$$ \Tr \lf| \ff{z-D_n^\eps} \ri|^2 \leq \frac{n}{|\text{Im}(z)|^2}.$$

Hence, the right hand side of \eqref{TracesProduct} is bounded by 

\begin{align*}
& \f{n^{\frac12}}{|\text{Im}(z)|} \ \E\lf[ \sum_{l=1}^n \lf(\ff{z-\wtD_n^\eps}\wtX_n\ri)^{3}_{i,l} \lf(\ovl{\ff{z-\wtD_n^\eps}\wtX_n}\ri)^{3}_{i,l} \ri]^\frac12
\\
&\leq \f{n^{\frac12}}{|\text{Im}(z)|} \ \E\left[  \sum_{j,k,l,m,p=1}^n \frac{ (\wtX_n)_{i,j} \ (\wtX_n)_{j,k} \ (\wtX_n)_{k,l} \ \ovl{(\wtX_n)_{i,m}} \ \ovl{(\wtX_n)_{m,p}} \ \ovl{(\wtX_n)_{p,l}}}{\lf|z-\widetilde{\lam}_n^\eps(i)\ri|^2 \ \lf(z-\widetilde{\lam}_n^\eps(j)\ri) \ \lf(z-\widetilde{\lam}_n^\eps(k)\ri) \ \lf(\ovl{z-\widetilde{\lam}_n^\eps(m)}\ri) \ \lf(\ovl{z-\widetilde{\lam}_n^\eps(p)}\ri) }\right]^{\ff{2}}.
\end{align*}


Recall that the diagonal of the matrix $\wtX_n$ is null, hence the denominators of the terms of the previous sum are independent from the numerators. Moreover the expectation of the numerators are null except when
the set of indices $\{(i,j), (j,k), (k,l)\}$ are equal to the set $\{(i,m), (m,p), (p,l)\}$. Therefore, the complexity of the previous sum is $O(n^3)$. 

Moreover, for all indices $j,k,l,m,p$, 
$$\E\lf[ \frac1{\lf|z-\widetilde{\lam}_n^\eps(i)\ri|^2 \ \lf(z-\widetilde{\lam}_n^\eps(j)\ri) \ \lf(z-\widetilde{\lam}_n^\eps(k)\ri) \ \lf(\ovl{z-\widetilde{\lam}_n^\eps(m)}\ri) \ \lf(\ovl{z-\widetilde{\lam}_n^\eps(p)}\ri) } \ri]  \leq \frac{1}{|\text{Im}(z)|^6}$$
and since the $L^{6}$ norm of the entries of $\sqrt nX$ is finite, we get that, uniformly in the indices $j,k,l,m,p$,
$$
\E\left[  (\wtX_n)_{i,j} \ (\wtX_n)_{j,k} \ (\wtX_n)_{k,l} \ \ovl{(\wtX_n)_{i,m}} \ \ovl{(\wtX_n)_{m,p}} \ \ovl{(\wtX_n)_{p,l}}\right] = O(n^{-3}).
$$
Hence,
$$ 
\E\left[  \sum_{j,k,l,m,p=1}^n \frac{ (\wtX_n)_{i,j} \ (\wtX_n)_{j,k} \ (\wtX_n)_{k,l} \ \ovl{(\wtX_n)_{i,m}} \ \ovl{(\wtX_n)_{m,p}} \ \ovl{(\wtX_n)_{p,l}}}{\lf|z-\widetilde{\lam}_n^\eps(i)\ri|^2 \ \lf(z-\widetilde{\lam}_n^\eps(j)\ri) \ \lf(z-\widetilde{\lam}_n^\eps(k)\ri) \ \lf(\ovl{z-\widetilde{\lam}_n^\eps(m)}\ri) \ \lf(\ovl{z-\widetilde{\lam}_n^\eps(p)}\ri) }\right]^{\ff{2}} \leq \frac{C}{|\text{Im}(z)|^3}.
$$

Therefore,
\[ \E\left[  \left((z-\wtD_n^\eps)^{-1} \wtX_n (z-\wtD_n^\eps)^{-1}\wtX_n (z-\wtD_n^\eps)^{-1}\wtX_n (z- D_n^\eps)^{-1} \right)^2_{ii} \right]^{\ff{2}} \leq  \f{C n^\frac12}{|\text{Im}(z)|^4}.\]
\epr

Gathering Formulas \eqref{Pi_nInterest}, \eqref{diagonalExpansion} and Claims $1$, $2$ and $3$, we prove Lemma \ref{lemmaPi_n}.
    
\subsection{From Stieltjes transform to $\cC^7$ functions}\label{subsec1312}
 
Now, let  $\vfi$ be  a   $\mc{C}^{7}$ function on $\R$ with bounded seventh derivative and let us introduce the almost analytic extension of degree $7$ of $\vfi$ defined by \begin{equation*}
\forall z=x+\ii y\in\C, \qquad \widetilde \vfi_6(z) \ := \ \sum_{k = 0}^6\frac{1}{k !} (\ii y)^k \vfi^{(k)}(x)\,.
\end{equation*}
An elementary computation gives, by successive cancellations, that
\be\label{eq23111614112} \bar \partial \widetilde \vfi_6(z) = \ff{2}\lf( \partial_x  + \ii \partial_y \ri) \widetilde \vfi_6(x+\ii y) = \ff{2\times 6!} (\ii y)^6 \vfi^{(7)}(x).\ee

Furthermore, by Helffer-Sj\"ostrand formula \cite[Propo. 9]{BGEM}, for $\chi \in \mathcal{C}^\infty_c(\C;[0,1])$   a smooth cutoff function with value one on the support of $\vfi$,
\begin{equation*}
\vfi(\cdot) \;=\; -\frac{1}{\pi} \int_{\C} \frac{\bar \partial (\widetilde \vfi_6(z) \chi(z))}{y^6}y^6\varphi_z(\cdot) \, \ud^2 z\,
\end{equation*}
where $\ud^2 z$ denotes the Lebesgue measure on $\C$.

Note that by \eqref{eq23111614112}, $z \mapsto \one_{y\ne 0} \frac{\bar \partial (\widetilde \vfi_6(z) \chi(z))}{y^6}$ is a continuous compactly supported function and that $z\in \C\mapsto  \one_{y\ne 0} y^6\vfi_z\in \mc{C}^1_{\op{b}}$ is continuous, hence,
\[ \Pi_n(\vfi) = \frac{1}{\pi} \int_{\C}  \f{\bar \partial (\widetilde \vfi_6(z) \chi(z))}{y^6} \ y^6\Pi_n(\varphi_z) \, \ud^2 z. \]
Therefore, using the Cauchy-Schwarz inequality and the fact that $\chi$ has compact support at the second step, for a certain constant   $C$, we have
\begin{align*}
\E\lf( \lf| \Pi_n(\vfi) \ri|^2 \ri) &= \E\lf( \lf| \frac{1}{\pi} \int_{\C}  \f{\bar \partial (\widetilde \vfi_6(z) \chi(z))}{y^6} \ y^6\Pi_n(\varphi_z) \, \ud^2 z \ri|^2  \ri)
\\
& \leq C\E\lf( \int_{\C} \lf| \f{\bar \partial (\widetilde \vfi_6(z) \chi(z))}{y^6} \ y^6\Pi_n(\varphi_z) \ri|^2 \, \ud^2 z  \ri)
\\
&=C \int_{\C} \lf| \f{\bar \partial (\widetilde \vfi_6(z) \chi(z))}{y^6}\ri|^2 \ y^{12} \ \E\lf(\lf| \Pi_n(\varphi_z)\ri|^2\ri)  \, \ud^2 z\, .
\end{align*}
By \eqre{eq23111614112},  the function $\lf| \f{\bar \partial (\widetilde \vfi_6(z) \chi(z))}{y^6}\ri|^2$ is continuous and compactly supported and bounded by  $ C \|\vfi^{(7)}\|_\infty^2$ for some constant $C$. Besides,  by Lemma \ref{lemmaPi_n},  uniformly in $z$, 
\[ y^{12} \ \E\lf(\lf| \Pi_n(\varphi_z)\ri|^2\ri)  = O \lf( (1+y^6) \left( \eta_n^2 + \frac{1}{n}  + n\eps^2 \ri)  \ri). \]

We deduce that 
\[ \E\lf( \lf| \Pi_n(\vfi) \ri|^2 \ri) \leq C \int_{\C} \lf| \f{\bar \partial (\widetilde \vfi_6(z) \chi(z))}{y^6}\ri|^2 \ y^{12} \ \E\lf(\lf| \Pi_n(\varphi_z)\ri|^2\ri)  \, \ud^2 z\, = \, O\left(\|\vfi^{(7)}\|_\infty^2 \ \left(\eta_n^2 + \frac{1}{n} + n\eps^2\right)\right) ,\] which closes the proof of \theo{theoMain0917}.

\section{Proof of \theo{theoMain0918}}\label{sectionHeuristicproof}

Let us start with the study of the term $\op{Card} \{ j\,:\, |\lam_j^\eps - f(x)| < \al_n  \}$. By Weyl's inequalities on the eigenvalues of sum of operators (see \cite[Cor. 4.3.15.]{horn1990matrix}),  the ordered  eigenvalues of the $D_n^{\eps}$ and $D_n$ do not differ by more than $\eps  \|X_n\|_{\text{op}}$. Therefore, we have, with \pro tending to one,   \begin{align*}&\op{Card}\left\{ j\,:\, |\lam_j - f(x)| < \al_n - \eps  \|X_n\|_{\text{op}} \right\}\le \\ &\qquad\qquad \qquad\qquad \qquad\qquad \op{Card}\left\{ j\,:\, |\lam_j^\eps - f(x)| < \al_n \right\}\le \\
&\qquad\qquad  \qquad\qquad\qquad\qquad \qquad\qquad \qquad\qquad\op{Card}\left\{ j\,:\, |\lam_j - f(x)| < \al_n + \eps  \|X_n\|_{\text{op}} \right\}.\end{align*}

Since the variances of $x_{i,j}^n$ are uniformly bounded by $\frac{C}{n}$ for some $C>0$, and that the entries $x^n_{i,j}$ have finite moment of order $6$, the assumptions of Lemma \ref{LemOpNorm} of the Appendix are satisfied by Tchebychev inequality and there exists a constant $C$ such that $\mathbb{P}\left( \|X_n\|_{\text{op}} > C\right)$ converges to zero.
Hence, since $\eps \ll \alpha_n$, the cardinality of  $\left\{ j\,:\, |\lam_j - f(x)| < \al_n - \eps \|X_n\|_{\text{op}} \right\}$ and of $\left\{ j\,:\, |\lam_j - f(x)| < \al_n + \eps  \|X_n\|_{\text{op}} \right\}$ are asymptotically equal to $2n\al_n\rho(x)(1+o(1))$.

Henceforth, for any measure $\mu$ and integrable function $\vfi$, we use the convenient notation, $\mu(\vfi) := \int \vfi \ud \mu$.

Let us turn to the estimation of the sum $\sum_{\left\{ j\,:\, |\lam_j^\eps - f(x)| < \al_n \right\}}  |\lan\bu_{j}^\eps,\bee_{\lfloor nx_0 \rfloor}\ran|^2$. 
Denoting $\vfi_{x, \al_n}(t) := \one_{t\in [ f(x)-\al_n, \  f(x)+\al_n]}$, the previous sum is nothing but $\mu_{n,\bee_{\lfloor n x_0\rfloor}}^\eps(\vfi_{x, \al_n})$. We want to apply Theorem \ref{theoMain0917}, but since $\vfi_{x, \al_n}$ is not smooth, we bound it from above and below after introducing some $\omega_n \ll \al_n$ we will calibrate further.

With the use of a decreasing  smooth function satisfying $\psi_{|\mathbb{R}_{-}} = 1$ and $\psi_{|[1, \infty)} = 0$, we can bound $\vfi_{x, \al_n}$ by two smooth functions $\vfi_{x, \al_n, \omega_n}^-$ and $\vfi_{x, \al_n, \omega_n}^+$ defined by 
\beq \vfi_{x, \al_n, \omega_n}^-(t) &:=& \psi\left( 1 + \frac{t-f(x) - \al_n}{\omega_n} \right) \psi\left( 1 - \frac{ t-f(x) + \al_n}{\omega_n} \right),
\\
\vfi_{x, \al_n, \omega_n}^+(t) &:=& \psi\left(\frac{t-f(x) - \al_n}{\omega_n} \right) \psi\left(- \frac{ t-f(x) + \al_n}{\omega_n} \right).\eeq

The properties of these functions are illustrated on the picture below.

\begin{center}
~\hspace{1.3cm}\includegraphics[width=0.7\linewidth]{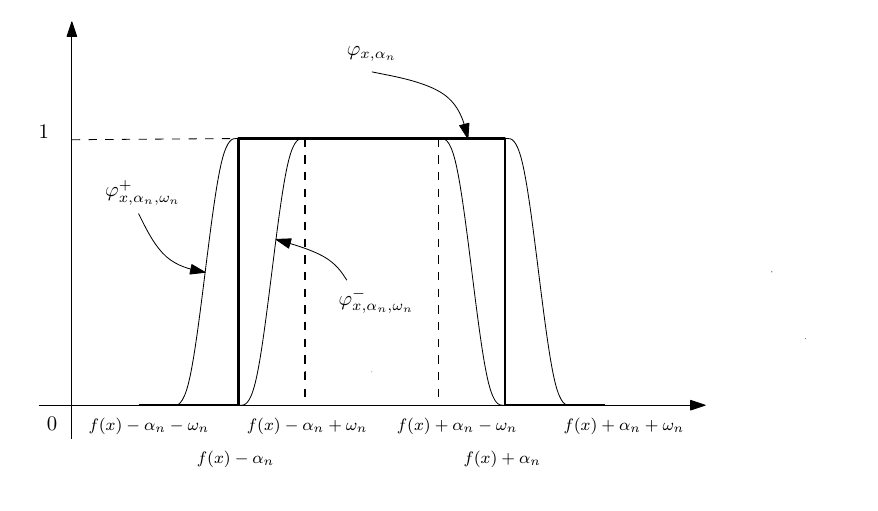}
\end{center}

Since $\mu_{n,\bee_{\lfloor n x_0\rfloor}}^\eps$ is a positive measure our quantity of interest $\mu_{n,\bee_{\lfloor n x_0\rfloor}}^\eps(\vfi_{x, \al_n})$ is bounded respectively from below and above by $\mu_{n,\bee_{\lfloor n x_0\rfloor}}^\eps(\vfi^-_{x, \al_n,\omega_n})$ and by $\mu_{n,\bee_{\lfloor n x_0\rfloor}}^\eps(\vfi^+_{x, \al_n, \omega_n})$. Therefore, we just have to prove that each of them is asymptotically equal in probability to $$2\eps^2\al_n \frac{\tau(f(x_0), f(x))}{\left( f(x) - f(x_0) \right)^2} \rho(f(x)).$$ We examine all the quantities of Theorem \ref{theoMain0917}. Obviously, since $\omega_n \ll \al_n$ and since the support of $\vfi_{x, \al_n, \omega_n}^-$  and $ \vfi_{x, \al_n, \omega_n}^+$ both avoid $f(x_0)$, the deterministic quantities $\eps^2 \Xi_{f(x_0)}(\vfi_{x, \al_n, \omega_n}^-)$ and $\eps^2 \Xi_{f(x_0)}(\vfi_{x, \al_n, \omega_n}^+)$ are asymptotically equal to the desired quantity announced before.

Now, since $\lam_{\lfloor nx_0 \rfloor}+\f{\eps}{\sqrt n} x_{\lfloor nx_0 \rfloor, \lfloor nx_0 \rfloor}$ converges in probability to $f(x_0)$ which is outside the support of $\vfi_{x, \al_n, \omega_n}^-$ and $\vfi_{x, \al_n, \omega_n}^+$ the quantity $\frac{n\eps^{-2}}{\al_n} \vfi_{x, \al_n, \omega_n}^{\pm}(\lam_{\lfloor nx_0 \rfloor}+\f{\eps}{\sqrt n} x_{\lfloor nx_0 \rfloor, \lfloor nx_0 \rfloor})$ converges in probability to zero.

Finally, the error term in Theorem \ref{theoMain0917} is, in $L^2$, of order $\eps^2 \| \vfi_{x, \al_n, \omega_n}^{\pm (7)} \|_\infty  (n^{\frac12} \eps+n^{-\frac12}+\eta_n)$ which in turn is of order $ \frac{\eps^2}{\omega_n^7} \left( n^{\frac12}\eps + n^{-\frac12} + \eta_n \right) $. It is now time to calibrate $\omega_n$ such that $\omega_n \ll \al_n$ and $\frac{\eps^2}{\omega_n^7} \left( n^{\frac12}\eps + n^{-\frac12}+ \eta_n \right) \ll \eps^2 \al_n$. This is possible if and only if $\al_n \gg \max\{(n^\frac{1}{2} \eps)^{\frac18}, n^{-\frac{1}{16}}, \eta_n^{\frac18}\}$.
This  closes the proof of \theo{theoMain0918}.

\section{Appendix}

This appendix is devoted to the control of the operator norm of the random matrix $X_n$ that we use in the proof of Theorem \ref{theoMain0918}. We did not find any reference for the following lemma in the literature, so we give a proof.   
  \blem\label{LemOpNorm} Let $H=(H_{ij})_{1\le i,j\le n}$ be an $n\ti n$ random Hermitian matrix satisfying:
  \bgt\ite The random variables $(H_{ij})_{1\le i\le j\le n}$ are independent, centered and satisfy $$n\E |H_{ij}|^2\le 1.$$
  \ite For some constants $C_0>0$, $\al>4$, we have, for any $t>0$, $$\p(\sqrt{n} |H_{ij}|\ge t)\le C_0 t^{-\al}.$$
  \ent
  Then  for any $\eps>0$, there is $C$ depending only on $\eps$,  $C_0$ and  $\al$ \st $$\p(\|H\|\ge 2+\eps)\le Cn^{-\f{\al-4}{4}}.$$
\elem
 
 \bpr Let $H':=(H_{ij}\one_{|H_{ij}|< n^{-\bet}})_{1\le i,j\le n}$ for $\bet:=\f{\al-4}{4\al}$.
 Note that by the union bound, $$\p(H\ne H')\le n^2C_0n^{-\al(1/2-\bet)}=C_0n^{-\f{\al-4}{4}},$$ so that it is enough to prove the result for $H'$ instead of $H$. The matrix $H'$ satisfies the assumptions of \cite[Th. 2.6]{FBGCBAK2016} for $\ka=1$ and $q:=n^{\max\{1/13, \bet\}}$, so, by this theorem and \cite[Eq. (2.4)]{FBGCBAK2016}, we know that for some universal positive constants $C_1,c$, for any $t\ge 0$, $$\p(\|H'\|\ge 2+C_1\sqrt{\log n}/q+t)\;\le\; 2\me^{-cq^2t^2}, 
 $$ which allows to conclude.
 \epr

\noindent{\bf Acknowledgments.} The authors wish to thank the  anonymous referees for their careful reading, their useful advices and for pointing out a mistake in the proof of Lemma \ref{lemmaPi_n}.

\end{document}